\documentclass[10 pt, draft]{article}
\title{The complex of partial bases for $F_n$ and finite generation of the Torelli subgroup of $\Aut(F_n)$}
\author{Matthew Day\footnote{Supported in part by an NSF postdoctoral fellowship}\ \ and Andrew Putman\footnote{Supported in part by NSF grant DMS-1005318}}
\usepackage{amsmath}
\usepackage{amssymb}
\usepackage{amsthm}
\usepackage{epsfig}
\usepackage{geometry}
\usepackage{amsfonts}
\usepackage{calc}
\usepackage{amscd}
\usepackage[font=small]{caption}
\usepackage{MnSymbol}
\usepackage{booktabs}
\usepackage{url}

\theoremstyle{plain}
\newtheorem{theorem}{Theorem}[section]
\newtheorem{maintheorem}{Theorem}

\newtheorem{lemma}[theorem]{Lemma}
\newtheorem*{unnumberedlemma}{Lemma}

\newtheorem{corollary}[theorem]{Corollary}
\newtheorem{conjecture}[theorem]{Conjecture}

\newtheorem{claim}{Claim}
\newtheorem{step}{Step}

\newcommand\BeginClaims{\setcounter{claim}{0}}

\newcommand\BeginSteps{\setcounter{step}{0}}

\theoremstyle{definition}
\newtheorem*{definition}{Definition}

\newtheorem*{outline}{Outline and conventions}

\theoremstyle{remark}
\newtheorem*{remark}{Remark}

\newtheorem*{warning}{Warning}

\DeclareMathOperator{\Hom}{Hom}

\DeclareMathOperator{\Ker}{ker}

\DeclareMathOperator{\Mod}{Mod}

\DeclareMathOperator{\IA}{IA}

\DeclareMathOperator{\GL}{GL}

\newcommand\Curves{\ensuremath{\mathcal{C}}}
\newcommand\CNosep{\ensuremath{\mathcal{C}^{\text{nosep}}}}

\newcommand\Z{\ensuremath{\mathbb{Z}}}

\DeclareMathOperator{\HH}{H}

\DeclareMathOperator{\Cayley}{Cay}
\DeclareMathOperator{\Star}{star}
\DeclareMathOperator{\Link}{link}
\DeclareMathOperator{\Max}{max}

\DeclareMathOperator{\Aut}{Aut}
\DeclareMathOperator{\Out}{Out}

\newcommand\Span[1]{\ensuremath{\langle #1 \rangle}}

\DeclareMathOperator{\Dim}{dim}
\newcommand\Set[2]{\ensuremath{\{\text{#1 $|$ #2}\}}}

\newcommand\Mul[2]{\ensuremath{\text{M}_{#1,#2}}}
\newcommand\Mulcomm[3]{\ensuremath{\text{M}_{#1,[#2,#3]}}}
\newcommand\Con[2]{\ensuremath{\text{C}_{#1,#2}}}
\DeclareMathOperator{\Rank}{Rank}
\newcommand\Conj[1]{\ensuremath{\lsem #1 \rsem}}

\DeclareMathOperator{\SAut}{SAut}
\newcommand\NormalSpan[1]{\ensuremath{\llangle #1 \rrangle}}
\newcommand\PartialBases[1]{\ensuremath{\mathcal{B}_{#1}}}
\newcommand\PartialBasesZ[1]{\ensuremath{\mathcal{B}_{#1}(\Z)}}
\newcommand\PartialBasesZZ[2]{\ensuremath{\mathcal{B}_{#1}^{#2}(\Z)}}
\newcommand\BKer[2]{\ensuremath{\mathcal{K}_{#1}(#2)}}
\DeclareMathOperator{\BKerr}{\mathcal{K}}

\begin{document}

\maketitle

\begin{abstract}
We study the complex of partial bases of a free group, which is
an analogue for $\Aut(F_n)$ of the curve complex for the mapping class
group.  We prove that it is connected and simply connected, and we
also prove that its quotient by the Torelli subgroup of $\Aut(F_n)$
is highly connected.  Using these results, we give a new, topological
proof of a theorem of Magnus that asserts that the Torelli subgroup
of $\Aut(F_n)$ is finitely generated.
\end{abstract}

\section{Introduction}
Let $\Sigma_g$ be a compact orientable genus $g$ surface and let $\Mod_g$ be its mapping class group.
One of the most important and ubiquitous objects associated to $\Mod_g$ is the
{\it curve complex} $\Curves_g$.  By definition, this is the simplicial complex whose 
simplices are sets of homotopy classes of non-nullhomotopic simple closed curves on 
$\Sigma_g$ that can be realized disjointly.  It can be viewed as an analogue for $\Mod_g$ of
the Tits building of an algebraic group.  The space $\Curves_g$ has many
remarkable properties; for instance, Harer \cite{HarerStability} showed that $\Curves(\Sigma_g)$ is 
homotopy equivalent to a bouquet of spheres and Masur-Minsky \cite{MasurMinskyHyperbolic}
showed that $\Curves(\Sigma_g)$ is $\delta$-hyperbolic.

There is a useful analogy between the mapping class group of a surface and the automorphism group of
a free group $F_n$ on $n$ letters.  Because of this, there have been many proposals for analogues of
the curve complex for $\Aut(F_n)$ (for instance, see 
\cite{BestvinaFeighn, HatcherSphereComplex, HatcherVogtmannFreeFactors, HatcherVogtmannCerf, KapovichLustig}).
The purpose of this paper is to prove some topological results about one of these proposed
complexes (the {\em complex of partial bases} $\PartialBases{n}$; see below).  We apply these results to give
a quick proof of a classical theorem of Magnus which provides generators for the {\em Torelli subgroup} 
$\IA_n < \Aut(F_n)$, which is the kernel of the natural homomorphism from $\Aut(F_n)$ to 
$\Aut(F_n^{\text{ab}}) \cong \Aut(\Z^n) \cong \GL_n(\Z)$.

Our complex is inspired by a subcomplex of $\Curves_g$, the {\em nonseparating curve complex}.
This is the subcomplex $\CNosep_g$ of $\Curves_g$ consisting of simplices $\{\gamma_1,\ldots,\gamma_k\}$ 
such that the $\gamma_i$ can be realized by simple closed curves whose union does not
separate $\Sigma_g$.  The complex $\CNosep_g$ was introduced by Harer \cite{HarerStability} and 
plays an important role in both homological stability results for $\Mod_g$ and its subgroups (see 
\cite{HarerStability, PutmanH2}) and the second author's approach to the Torelli subgroup
of $\Mod_g$ (see \cite{PutmanCutPaste, PutmanInfinite, PutmanJohnson, PutmanAbel}).

To motivate the definition of our complex, we start by giving an algebraic characterization
of $\CNosep_g$.  There is a bijection between free homotopy classes of oriented closed curves on $\Sigma_g$ and
conjugacy classes in $\pi_1(\Sigma)$.  We will call a generating set 
$\{\alpha_1,\beta_1,\ldots,\alpha_g,\beta_g\}$ for $\pi_1(\Sigma_g)$ a {\em standard basis}
if it satisfies the surface relation $[\alpha_1,\beta_1] \cdots [\alpha_g,\beta_g] = 1$.  We then
have the following folklore result, whose proof is an exercise in using the fact that $\Mod_g$ is
the outer automorphism group of $\pi_1(\Sigma_g)$.  If $G$ is a group and $g \in G$, then denote
by $\Conj{g}$ the conjugacy class of $g$.

\begin{unnumberedlemma}
A set $\{c_1,\ldots,c_n\}$ of conjugacy classes in $\pi_1(\Sigma_g)$ corresponds to a simplex
of $\CNosep_g$ (with some orientation on each curve in the simplex) if and only if there
exists a standard basis $\{\alpha_1,\beta_1,\ldots,\alpha_g,\beta_g\}$ for $\pi_1(\Sigma_g)$ 
such that $c_i = \Conj{\alpha_i}$ for $1 \leq i \leq k$.
\end{unnumberedlemma}

\noindent
This suggests the following definition.

\begin{definition}
Fix $n \geq 1$.  A {\em partial basis} for $F_n$ consists of elements $\{v_1,\ldots,v_k\} \subset F_n$, with $1\leq k\leq n$, such that there exists $v_{k+1},\ldots,v_n \in F_n$
with $\{v_1,\ldots,v_n\}$ a free basis for $F_n$.
The {\em complex of partial bases} of $F_n$, denoted $\PartialBases{n}$, is the simplicial complex whose $(k-1)$-simplices are unordered sets $\{\Conj{v_1},\ldots,\Conj{v_k}\}$ such that $\{v_1,\ldots,v_k\} \subset F_n$ is a partial basis for $F_n$.
\end{definition}

We consider any complete basis for $F_n$ to be a partial basis, but by convention the empty set is not a partial basis.

\begin{remark}
Our theorems are about topological properties of $\PartialBases{n}$ and algebraic properties of $\Aut(F_n)$.
However, the geometry of $\PartialBases{n}$ and related spaces has recently received much 
attention in the literature.  The primitivity graph of Kapovich--Lustig \cite{KapovichLustig} is 
the $1$-skeleton of $\PartialBases{n}$; that paper discusses the geometric relationships 
between the primitivity graph (and therefore $\PartialBases{n}$) and several other curve-complex analogues. 
Kapovich--Lustig also prove that $\PartialBases{n}$ has infinite diameter for $n \geq 2$.  Very recently
Bestvina--Feighn proved that the {\em free factor complex} is $\delta$-hyperbolic.  This
complex is quasi-isometric to $\PartialBases{n}$ for $n \geq 3$, so this implies that $\PartialBases{n}$
is $\delta$-hyperbolic as well.  Other recent papers on the geometry of analogues of the curve
complex include \cite{HandelMosher, SabalkaSavchuk}.
\end{remark}

Our two main results about $\PartialBases{n}$ are as follows.

\begin{maintheorem}
\label{theorem:basesconnected}
The space $\PartialBases{n}$ is connected for $n \geq 2$ and $1$-connected for $n \geq 3$.
\end{maintheorem}

\begin{maintheorem}
\label{theorem:basesiaconnected}
The space $\PartialBases{n} / \IA_n$ is $(n-2)$-connected.
\end{maintheorem}

\noindent
The analogues of Theorems \ref{theorem:basesconnected} and \ref{theorem:basesiaconnected} for $\CNosep_g$ are due to Harer \cite{HarerStability}
and the second author \cite{PutmanInfinite}, respectively.  In fact, Harer proved $\CNosep_g$ is $(g-2)$-connected, which leads us to make
the following conjecture.

\begin{conjecture}
The space $\PartialBases{n}$ is $(n-2)$-connected.
\end{conjecture}

\begin{remark}
In their papers \cite{HatcherVogtmannFreeFactors, HatcherVogtmannCerf}, Hatcher-Vogtmann prove that various simplicial complexes
built out of free factors of $F_n$ are highly connected.  Their proofs are quite different from our proof of Theorem \ref{theorem:basesconnected}
and do not seem to apply (at least directly) to $\PartialBases{n}$.  Also, their complexes are not $1$-connected for $n=3$, which renders
them unsuitable for our application below to $\IA_n$ (the induction needs the case $n=3$ to get started).
\end{remark}

Next, we need another definition.  

\begin{definition}
Let $\{v_1,\ldots,v_n\}$ be a fixed free basis for $F_n$.  Choose $1 \leq i \leq n$, and let 
$w \in F_n$ be an element of
the subgroup of $F_n$ spanned by $\Set{$v_j$}{$j \neq i$}$.
\begin{itemize}
\item For $e \in \{1,-1\}$, let $\Mul{v_i^e}{w} \in \Aut(F_n)$ denote the automorphism that takes $v_i^e$ to $wv_i^e$ and fixes $v_j$ for
$j \neq i$.
\item Let $\Con{v_i}{w} \in \Aut(F_n)$ denote the automorphism that takes $v_i$ to $w v_i w^{-1}$ and fixes $v_j$ for $j \neq i$.
\end{itemize}
\end{definition}

\begin{remark}
Observe that $\Con{v_i}{w} = \Mul{v_i}{w} \Mul{v_i^{-1}}{w}$.
\end{remark}

\begin{remark}
In the literature, it is common to instead work with the elements $\Mul{v_i^e}{w}' \in \Aut(F_n)$ that take $v_i^e$ to $v_i^e w$ and fix
$v_j$ for $j \neq i$.  If one used this convention, then it would make sense to define $\Con{v_i}{w}$ to take $v_i$ to $w^{-1} v_i w$; however,
this conjugation convention would make conjugation a right action rather than a left action.  We prefer to work with left actions.
\end{remark}

\noindent
Using the complex $\PartialBases{n}$ and Theorems \ref{theorem:basesconnected} and \ref{theorem:basesiaconnected},
we give a new proof of the following theorem of Magnus.

\begin{maintheorem}[{Magnus, \cite{Magnus}}]
\label{theorem:magnusfinite}
Fix a free basis $\{v_1,\ldots,v_n\}$ for $F_n$.  The group $\IA_n$ is then generated by the finite set
\[\Set{$\Mul{v_i}{[v_j,v_k]}$}{$1 \leq i,j,k \leq n$, $i \neq j$, $i \neq k$, $j \neq k$} \cup \Set{$\Con{v_i}{v_j}$}{$1 \leq i,j \leq n$, $i \neq j$}.\]
\end{maintheorem}

\begin{remark}
The usual statement of Theorem \ref{theorem:magnusfinite} contains the elements $\Mul{v_i}{[v_j,v_k]}'$ 
instead of $\Mul{v_i}{[v_j,v_k]}$.  The two generating sets are equivalent due to the formula
$$\Mul{v_i}{[v_j,v_k]}' = \Con{v_i}{v_j}^{-1} \Con{v_i}{v_k}^{-1} \Con{v_i}{v_j} \Con{v_i}{v_k} \Mul{v_i}{[v_j,v_k]}.$$
\end{remark}

\begin{remark}
It was proven independently by Cohen-Pakianathan \cite{CohenPakianathan}, Farb \cite{Farb}, and Kawazumi \cite{Kawazumi} that
the size of the generating set in Theorem \ref{theorem:magnusfinite} is as small as possible.
\end{remark}

Magnus's original proof of Theorem \ref{theorem:magnusfinite} had two steps.  The first and 
most difficult is the following result.

\begin{maintheorem}[{Magnus, \cite{Magnus}}]
\label{theorem:magnus}
Fix a free basis $\{v_1,\ldots,v_n\}$ for $F_n$.  The group $\IA_n$ is then normally 
generated as a subgroup of $\Aut(F_n)$ by the single element $\Con{v_1}{v_2}$.
\end{maintheorem}

\noindent
In \S \ref{section:magnus} below, we give a 
topological proof of Theorem \ref{theorem:magnus} using the complex of partial bases.
Following Magnus, one can then derive Theorem \ref{theorem:magnusfinite} from Theorem \ref{theorem:magnus}
by showing that the subgroup of $\Aut(F_n)$ generated by the generating set in Theorem \ref{theorem:magnusfinite}
is normal.  For completeness, we give the details of this argument in Appendix \ref{appendix:magnus}.

\begin{remark}
There is another (quite different) topological proof of Theorem \ref{theorem:magnus} in a recent 
paper of Bestvina-Bux-Margalit \cite{BestvinaBuxMargalit}.
That paper also contains in an appendix a sketch of Magnus's 
derivation of Theorem \ref{theorem:magnusfinite} from Theorem \ref{theorem:magnus}.
\end{remark}

\begin{remark}
An analogue of Theorem \ref{theorem:magnus} for $\Mod_g$ was proven by Powell \cite{PowellTorelli}, 
following work of Birman \cite{BirmanSiegel}.
Their proof resembled Magnus's proof of Theorem \ref{theorem:magnus}, though the details were far more 
complicated.  Later, an analogue of 
Theorem \ref{theorem:magnusfinite} for $\Mod_g$ was proven by Johnson \cite{JohnsonFinite}, again 
following Magnus's derivation of
Theorem \ref{theorem:magnusfinite} from Theorem \ref{theorem:magnus} (though again the details are much more 
complicated).  Recently, the second
author \cite{PutmanCutPaste} gave a topological proof of Birman-Powell's theorem using the curve complex.  
The machinery of the current paper is set up so that the proof of Theorem \ref{theorem:magnus} follows the 
topological argument in \cite{PutmanCutPaste} closely.
\end{remark}

\begin{outline}
One of the key tools in this paper is an analogue of the Birman exact sequence for $\Aut(F_n)$ which the authors proved in \cite{DayPutmanBirman}.  This
is discussed in \S \ref{section:birman}.  Next, Theorem \ref{theorem:magnus} is proven in \S \ref{section:magnus} (assuming the truth of Theorems \ref{theorem:basesconnected}
and \ref{theorem:basesiaconnected}).  Finally, Theorems \ref{theorem:basesconnected} and \ref{theorem:basesiaconnected} are proven
in \S \ref{section:partialbases} and \S \ref{section:partialbasesia}.  The last three sections are largely independent of each other and can be read in any order.

Automorphisms act on the left and compose from right to left like functions.
Also, if $G$ is a group and $g,h \in G$, then we define $[g,h] = g h g^{-1} h^{-1}$.
\end{outline}

\section{The Birman exact sequence for $\Aut(F_n)$}
\label{section:birman}

A key tool in this paper is a type of Birman exact sequence for $\Aut(F_n)$ which the authors developed in \cite{DayPutmanBirman}.
This is an analogue of the classical Birman exact sequence for mapping class group (see \cite{BirmanSequence}).  
We begin with some notation.  For $v_1,\ldots,v_k \in F_n$, let
\[\Aut(F_n,\Conj{v_1},\ldots,\Conj{v_k}) = \Set{$\phi \in \Aut(F_n)$}{$\Conj{\phi(v_i)} = \Conj{v_i}$ for $1 \leq i \leq k$}.\]
Also, if $G$ is a group and $S \subset G$, then let $\Span{S} \subset G$ denote the subgroup
of $G$ generated by $S$ and $\NormalSpan{S} \subset G$ denote the normal subgroup of $G$
generated by $S$.

Let $\{v_1,\ldots,v_n\}$ be a basis for $F_n$.  For some $1 \leq k \leq n$, let 
$V = \NormalSpan{v_1,\ldots,v_k}$.  Since the group
$\Aut(F_n,\Conj{v_1},\ldots,\Conj{v_k})$ preserves $V$, we get an induced map
\[\Aut(F_n,\Conj{v_1},\ldots,\Conj{v_k}) \rightarrow \Aut(F_n/V)\]
which is clearly a split surjection.  Define $\BKer{n}{\Conj{v_1},\ldots,\Conj{v_k}}$ to be
its kernel, so we have a split short exact sequence
\begin{equation}
\label{eqn:birmanaut}
1 \longrightarrow \BKer{n}{\Conj{v_1},\ldots,\Conj{v_k}} \longrightarrow \Aut(F_n,\Conj{v_1},\ldots,\Conj{v_k}) \longrightarrow \Aut(F_n/V) \longrightarrow 1.
\end{equation}
The paper \cite{DayPutmanBirman} contains numerous results about $\BKer{n}{\Conj{v_1},\ldots,\Conj{v_k}}$; for instance, it is proven that
it is finitely generated but not finitely presentable, its abelianization is calculated, and an explicit infinite presentation
for it is constructed.  We will only need the following result.

\begin{theorem}[{\cite{DayPutmanBirman}}]
\label{theorem:birmankergen}
If $\{v_1,\ldots,v_n\}$ is a basis for $F_n$ and $1 \leq k \leq n$, then $\BKer{n}{\Conj{v_1},\ldots,\Conj{v_k}}$ is
generated by
\[\Set{$\Mul{v_i}{v_j}$}{$k+1 \leq i \leq n$, $1 \leq j \leq k$} \cup
\Set{$\Con{v_i}{v_j}$}{$1 \leq i \leq k$, $1 \leq j \leq n$, $i \neq j$}.\]
\end{theorem}

The exact sequence \eqref{eqn:birmanaut} is related to the following exact sequence for $\Aut(\Z^n)$.  Let $\overline{V} \subset \Z^n$
be a direct summand and define
\[\Aut(\Z^n,\overline{V}) = \Set{$\phi \in \Aut(\Z^n)$}{$\phi|_{\overline{V}} = 1$}.\]
For an appropriate choice of basis, $\Aut(\Z^n,\overline{V})$ consists of matrices in $\GL_n(\Z)$ with an
identity block in the upper left hand corner and a block of zeros in the lower left
hand corner.  We then have a split short exact sequence
\begin{equation}
\label{eqn:birmangl}
1 \longrightarrow \Hom(\Z^n/\overline{V},\overline{V}) \longrightarrow \Aut(\Z^n,\overline{V}) \longrightarrow \Aut(\Z^n/\overline{V}) \longrightarrow 1.
\end{equation}
Letting $\pi \colon \Z^n \rightarrow \Z^n/V$ be the projection, an element $\varphi \in \Hom(\Z^n/\overline{V},\overline{V})$
corresponds to the element of $\Aut(\Z^n,\overline{V})$ that takes $x \in \Z^n$ to $x + \varphi(\pi(x))$.

The exact sequences \eqref{eqn:birmanaut} and \eqref{eqn:birmangl} are related by the following lemma.

\begin{lemma}[Images of stabilizers]
\label{lemma:autstabimage}
Fix bases $\{v_1,\ldots,v_n\}$ and $\{\overline{v}_1,\ldots,\overline{v}_n\}$ for $F_n$ and $\Z^n$, 
respectively, such that
$\pi(v_i) = \overline{v}_i$ for $1 \leq i \leq n$.  Choose some $1 \leq k \leq n$.  Set 
$\overline{V} = \Span{\overline{v}_1,\ldots,\overline{v}_k} \subset \Z^n$ and 
$V = \NormalSpan{v_1,\ldots,v_k} \subset F_n$.
There is then a commutative diagram of split short exact sequences
\[\minCDarrowwidth20pt\begin{CD}
1 @>>> \BKer{n}{\Conj{v_1},\ldots,\Conj{v_k}} @>>> \Aut(F_n,\Conj{v_1},\ldots,\Conj{v_k})  @>>> \Aut(F_n/V)             @>>> 1\\
@.     @VVV                                        @VVV                                         @VVV                         @.\\
1 @>>> \Hom(\Z^n/\overline{V},\overline{V})   @>>> \Aut(\Z^n,\overline{V})                 @>>> \Aut(\Z^n/\overline{V}) @>>> 1
\end{CD}\]
whose vertical maps are all surjective.
\end{lemma}

\begin{remark}
By a commutative diagram of split short exact sequences, we mean not only that the diagram commutes and that two sequences are split, but also
that the splitting is compatible with the commutative diagram in the obvious way.
\end{remark}

\begin{proof}[{Proof of Lemma \ref{lemma:autstabimage}}]
The only nonobvious claims are the surjectivity of the vertical maps.  The fact that the map
$\Aut(F_n/V) \rightarrow \Aut(\Z^n/\overline{V})$ is surjective is classical, while the fact that the map 
$\BKer{n}{\Conj{v_1},\ldots,\Conj{v_k}} \rightarrow \Hom(\Z^n/\overline{V},\overline{V})$
is surjective follows from the fact that the elements 
\[\Set{$\Mul{v_i}{v_j}$}{$k+1 \leq i \leq n$, $1 \leq j \leq k$} \subset \BKer{n}{\Conj{v_1},\ldots,\Conj{v_k}}\]
project to a basis for $\Hom(\Z^n/\overline{V},\overline{V})$.  The surjectivity of the map
$\Aut(F_n,\Conj{v_1},\ldots,\Conj{v_k}) \rightarrow \Aut(\Z^n,\overline{V})$ now follows from the five lemma.
\end{proof}

\noindent
Theorem \ref{theorem:birmankergen} and Lemma \ref{lemma:autstabimage} have the following corollary.

\begin{corollary}
\label{corollary:bkeriagen}
Let $\{v_1,\ldots,v_n\}$ be a basis for $F_n$.  Then for $1 \leq k \leq n$, the group
$\BKer{n}{\Conj{v_1},\ldots,\Conj{v_k}} \cap \IA_n$ is normally generated as a subgroup
of $\BKer{n}{\Conj{v_1},\ldots,\Conj{v_k}}$ by
\[\Set{$\Con{v_i}{v_j}$}{$1 \leq i \leq k$, $1 \leq j \leq n$, $i \neq j$}.\]
\end{corollary}
\begin{proof}
Let $\pi \colon F_n \rightarrow \Z^n$ be the abelianization map and let 
$\overline{V} = \Span{\pi(v_1),\ldots,\pi(v_k)}$.
Using Lemma \ref{lemma:autstabimage}, we have that $\BKer{n}{\Conj{v_1},\ldots,\Conj{v_k}} \cap \IA_n$
is the kernel of a surjective map
\[\Phi \colon \BKer{n}{\Conj{v_1},\ldots,\Conj{v_k}} \longrightarrow \Hom(\Z^n / \overline{V}, \overline{V}).\]
Let $S_{\BKerr}$ be the generating set for $\BKer{n}{\Conj{v_1},\ldots,\Conj{v_k}}$ given
by Theorem \ref{theorem:birmankergen} and let $S_{\BKerr,\IA} \subset S_{\BKerr}$ be the set we are trying 
to prove normally generates $\BKer{n}{\Conj{v_1},\ldots,\Conj{v_k}} \cap \IA_n$.  Also, let 
\[T = S_{\BKerr} \setminus S_{\BKerr,\IA} = \Set{$\Mul{v_i}{v_j}$}{$k+1 \leq i \leq n$, $1 \leq j \leq k$}.\]
Observe the following two facts.
\begin{itemize}
\item $S_{\BKerr,\IA} \subset \Ker(\Phi)$.
\item $\Phi|_{T}$ is injective and $\Phi(T)$ is a basis for the free abelian group 
$\Hom(\Z^n / \overline{V}, \overline{V})$.
\end{itemize}
From these two facts, it follows immediately that $\Ker(\Phi)$ is normally generated by the set
$S_{\BKerr,\IA} \cup \Set{$[t_1,t_2]$}{$t_1,t_2 \in T$}$.
We must show that the set $\Set{$[t_1,t_2]$}{$t_1,t_2 \in T$}$ is in the normal
closure of $S_{\BKerr,\IA}$.

Consider $\Mul{v_i}{v_j}, \Mul{v_{i'}}{v_{j'}} \in T$.  We want to express $[\Mul{v_i}{v_j}, \Mul{v_{i'}}{v_{j'}}]$
as a product of conjugates of elements of $S_{\BKerr,\IA}$.  If either $i \neq i'$ or $j = j'$, then
$[\Mul{v_i}{v_j}, \Mul{v_{i'}}{v_{j'}}] = 1$ and the claim is trivial (this uses the fact that
$\{v_i,v_{i'}\} \cap \{v_j, v_{j'}\} = \emptyset$).  Assume, therefore, that $i = i'$ and 
$j \neq j'$.  We then have the easily verified identity
\[[\Mul{v_i}{v_j}, \Mul{v_{i'}}{v_{j'}}] = \Mul{v_i}{[v_{j'}^{-1},v_{j}^{-1}]} = (\Mul{v_i}{v_j} \Con{v_i}{v_{j'}} \Mul{v_i}{v_j}^{-1}) \Con{v_i}{v_{j'}}^{-1},\]
and the corollary follows.
\end{proof}

\section{Generators for $\IA_n$}
\label{section:magnus}

We will now assume the truth of Theorems \ref{theorem:basesconnected} and \ref{theorem:basesiaconnected} and prove Theorem \ref{theorem:magnus}, which gives generators for $\IA_n$.  
The proof will closely follow the proof in \cite{PutmanCutPaste} of the analogous fact for
the mapping class group.  We start with a definition.

\begin{definition}
A group $G$ acts on a simplicial complex $X$ {\em without rotations} if for all
simplices $s$ of $X$, the stabilizer $G_s$ stabilizes $s$ pointwise.
\end{definition}

\noindent
For example, $\IA_n$ acts without rotations on $\PartialBases{n}$, but $\Aut(F_n)$ does not act
without rotations on $\PartialBases{n}$ if $n \geq 2$.

The key to our proof will be the following theorem of Armstrong \cite{ArmstrongGenerators}.  Our
formulation is a little different from Armstrong's; see
\cite[Theorem 2.1]{PutmanCutPaste} for a description of how to extract it from \cite{ArmstrongGenerators}.

\begin{theorem}[{Armstrong, \cite{ArmstrongGenerators}}]
\label{theorem:armstrong}
Let $G$ act without rotations on a $1$-connected simplicial complex $X$.  Then $G$ is generated by
the set
\[\bigcup_{v \in X^{(0)}} G_v\]
if and only if $X/G$ is $1$-connected.
\end{theorem}

\noindent
We now proceed to the proof.

\begin{proof}[Proof of Theorem \ref{theorem:magnus}]
Let $\{v_1,\ldots,v_n\}$ be a basis for $F_n$.  Our goal is to show that $\IA_n$ is normally
generated as a subgroup of $\Aut(F_n)$ by $\Con{v_1}{v_2}$.  We first observe that the conjugacy
class of $\Con{v_1}{v_2}$ is independent of the initial choice of basis.  It is thus enough to show
that $\IA_n$ is generated by the automorphisms that (for some choice of basis) conjugate one basis
element by another while fixing the rest of the basis.  We will call such an automorphism a
{\em basis-conjugating automorphism}.

The proof will be by induction on $n$.  The case $n=1$ is trivial, while the case $n=2$ follows
from a classical theorem of Nielsen \cite{NielsenPaper} that asserts that $\Out(F_2) \cong \GL_2(\Z)$, and
thus that $\IA_2$ consists entirely of inner automorphisms (see \cite[Proposition 4.5]{LyndonSchupp} for
a modern account of Nielsen's theorem).  Assume, therefore, that $n \geq 3$ and that the
result is true for all smaller $n$.  For $v \in F_n$, define
\[\IA_n(\Conj{v}) = \Set{$\phi \in \IA_n$}{$\Conj{\phi(v)} = \Conj{v}$}.\]
Using Theorems \ref{theorem:basesconnected} and \ref{theorem:basesiaconnected}, we can apply
Theorem \ref{theorem:armstrong} and deduce that $\IA_n$ is generated by the set
\[\bigcup_{\Conj{v} \in (\PartialBases{n})^{(0)}} \IA_n(\Conj{v}).\]
Consider some $\Conj{v} \in (\PartialBases{n})^{(0)}$.  Applying 
Lemma \ref{lemma:autstabimage}, we obtain a split short exact sequence
\[1 \longrightarrow \BKer{n}{\Conj{v}} \cap \IA_n \longrightarrow \IA_n(\Conj{v}) \longrightarrow \IA_{n-1} \longrightarrow 1.\]
The inductive hypothesis says that $\IA_{n-1}$ is generated by basis-conjugating automorphisms, and
Corollary \ref{corollary:bkeriagen} says that $\BKer{n}{\Conj{v}} \cap \IA_n$ is generated by
basis-conjugating automorphisms.  We conclude that $\IA_n(\Conj{v})$, and hence $\IA_n$, is generated
by basis-conjugating automorphisms, as desired.
\end{proof}

\section{The connectivity of $\PartialBases{n}$}
\label{section:partialbases}

In this section, we prove Theorem \ref{theorem:basesconnected}.  The proof uses a trick
that was introduced by the second author in \cite{PutmanConnectivity}.  We will need two results for
the proof.  For the first, let $\SAut(F_n)$ be the {\em special automorphism group} of $F_n$, i.e.\ the subgroup
of $\Aut(F_n)$ consisting of automorphisms whose images in $\Aut(\Z^n)$ have determinant $1$.  Also, if $S$ is a fixed basis
for $F_n$ and $a,b \in S^{\pm 1}$ satisfy $a \neq b^{-1}$, then denote by $w_{a,b}$ the automorphism of $F_n$ that
takes $a$ to $b^{-1}$, that takes $b$ to $a$, and fixes all the elements of $S^{\pm 1} \setminus \{a^{\pm 1}, b^{\pm 1}\}$.

\begin{theorem}[{Gersten, \cite{GerstenSpecial}}]
\label{theorem:gersten}
For $n \geq 3$, the group $\SAut(F_n)$ has the following presentation.  Let $S_{F}$ be a fixed basis for $F_n$.
\begin{itemize}
\item The generating set consists of
\[\Set{$\Mul{a}{b}$}{$a, b \in S_{F}^{\pm 1}$, $a \neq b^{\pm 1}$} \cup \Set{$\Con{a}{b}$}{$a,b \in S_{F}$, $a \neq b$} \cup \Set{$w_{a,b}$}{$a,b \in S_{F}^{\pm 1}$,$a \neq b^{\pm 1}$}.\]
\item The relations consist of the following.
\begin{enumerate}
\item $\Mul{a}{b} \Mul{a}{b^{-1}}=1$ for $a,b \in S_{F}^{\pm 1}$ with $a \neq b^{\pm 1}$.
\item $[\Mul{a}{b},\Mul{c}{d}] = 1$ for $a,b,c,d \in S_{F}^{\pm 1}$ with $b \neq a^{\pm 1}, c^{\pm 1}$ and $d \neq c^{\pm 1}, a^{\pm 1}$ and $a \neq c$.
\item $[\Mul{b}{a^{-1}},\Mul{c}{b^{-1}}] = \Mul{c}{a}$ for $a,b,c \in S_{F}^{\pm 1}$ with $a \neq b^{\pm 1},c^{\pm 1}$ and $b \neq c^{\pm 1}$.
\item $w_{a,b} = w_{a^{-1},b^{-1}}$ for $a,b \in S_{F}^{\pm 1}$ with $a \neq b^{\pm 1}$.
\item $w_{a,b} = \Mul{b^{-1}}{a^{-1}} \Mul{a^{-1}}{b} \Mul{b}{a}$ for $a,b \in S_{F}^{\pm 1}$ with $a \neq b^{\pm 1}$.
\item $w_{a,b}^4 = 1$ for $a,b \in S_{F}^{\pm 1}$ with $a \neq b^{\pm 1}$.
\item $\Con{a}{b} = \Mul{a}{b} \Mul{a^{-1}}{b}$ for $a,b \in S_{F}$ with $a \neq b$.
\end{enumerate}
\end{itemize}
\end{theorem}

\begin{remark}
The presentation above differs from the presentation in \cite{GerstenSpecial} in three ways.  First, our convention 
is to write functions on the left, while in \cite{GerstenSpecial} functions are written on the right.  
Second, our convention is that $\Mul{a}{b}$ multiplies $a$ on the left by $b$, while in \cite{GerstenSpecial} the
analogous elements multiply $a$ on the right by $b$.
Third, we have added the additional generators $w_{a,b}$ and $\Con{a}{b}$ together with relations 5 and 7, 
which express $w_{a,b}$ and $\Con{b}{a}$ in terms of the generators $\Mul{a}{b}$ of \cite{GerstenSpecial}.  
This will simplify our proof later on.
\end{remark}

The second result we will need is as follows.  For a partial basis $\{v_1,\ldots,v_k\}$ of $F_n$, define
\[\SAut(F_n,\Conj{v_1},\ldots,\Conj{v_k}) = \Set{$\phi \in \SAut(F_n)$}{$\Conj{\phi(v_i)} = \Conj{v_i}$ for $1 \leq i \leq k$}.\]

\begin{lemma}
\label{lemma:stabgen}
Let $S = \{v_1,\ldots,v_n\}$ be a fixed free basis of $F_n$.  Then for $1 \leq k \leq n$, the
group $\SAut(F_n,\Conj{v_1},\ldots,\Conj{v_k})$ is generated by the set
\[\Set{$\Mul{a}{b}$}{$a,b \in S^{\pm 1}$, $a \neq b^{\pm 1}$, and $a \notin \{v_1^{\pm 1},\ldots,v_k^{\pm 1}\}$} \cup
\Set{$\Con{a}{b}$}{$a,b \in S^{\pm 1}$ and $a \neq b^{\pm 1}$}\]
\end{lemma}
\begin{proof}
We will assume that $k<n$; the case of $k=n$ is similar but easier.
Let $T$ be the indicated generating set.  Also, let $I_n \in \Aut(F_n,\Conj{v_1},\ldots,\Conj{v_k})$ be the 
automorphism that takes $v_n$ to $v_n^{-1}$ and fixes $v_i$ for $1 \leq i < n$.  The short exact sequence
\[1 \longrightarrow \SAut(F_n,\Conj{v_1},\ldots,\Conj{v_k}) \longrightarrow \Aut(F_n,\Conj{v_1},\ldots,\Conj{v_k}) \longrightarrow \mathbb{Z}/2\mathbb{Z} \longrightarrow 1\]
splits via a splitting taking the generator of $\mathbb{Z}/2\mathbb{Z}$ to $I_n$.  In other words,
\[\Aut(F_n,\Conj{v_1},\ldots,\Conj{v_k}) = \SAut(F_n,\Conj{v_1},\ldots,\Conj{v_k}) \rtimes \Span{I_n}.\]
Observe that $I_n^{-1} s I_n \in T$ for all $s \in T$.  To prove the lemma, it therefore suffices
to show that $T \cup \{I_n\}$ generates $\Aut(F_n,\Conj{v_1},\ldots,\Conj{v_k})$

Let $G \subset \Aut(F_n,\Conj{v_1},\ldots,\Conj{v_k})$ be the subgroup consisting of 
automorphisms $\phi$ with the following properties.
\begin{itemize}
\item $\phi(v_i) = v_i$ for $1 \leq i \leq k$.
\item $\phi(v_i) \in \{v_{k+1}^{\pm 1}, \ldots, v_n^{\pm 1}\}$ for $k+1 \leq i \leq n$.
\end{itemize}
Making use of a deep theorem of McCool \cite{McCoolPresentation}, Jensen-Wahl \cite{JensenWahl} proved that 
$T \cup G$ generates $\Aut(F_n,\Conj{v_1},\ldots,\Conj{v_k})$.  It is an easy exercise 
to show that $G \cap \SAut(F_n,\Conj{v_1},\ldots,\Conj{v_k})$ is contained in the subgroup generated by 
$T$.  Since $G$ is generated by the union of $\{I_n\}$ and 
$G \cap \SAut(F_n,\Conj{v_1},\ldots,\Conj{v_k})$, we
conclude that $T \cup \{I_n\}$ generates $\Aut(F_n,\Conj{v_1},\ldots,\Conj{v_k})$,
as desired.
\end{proof}

\noindent
Finally, we recall the following definition.

\begin{definition}
Let $G$ be a group with a generating set $S$.  The {\em Cayley graph} of $G$, denoted
$\Cayley(G,S)$, is the graph with vertex set $G$ and with $g \in G$ connected by an edge to $g' \in G$
precisely when $g' = g s$ for some $s \in S^{\pm 1}$.
\end{definition}

\noindent
We can now prove Theorem \ref{theorem:basesconnected}.

\begin{proof}[{Proof of Theorem \ref{theorem:basesconnected}}]
Let $S_{F}=\{v_1,\ldots,v_n\}$ be a free basis for $F_n$ and let $S$ be the corresponding
generating set for $\SAut(F_n)$ from Theorem \ref{theorem:gersten}.  Observe that for $s \in S$,
either $\Conj{s(v_1)}=\Conj{v_1}$ or $\{v_1,s(v_1)\}$ forms a partial basis for $F_n$.  This implies
that the map
\begin{align*}
\SAut(F_n) &\rightarrow \PartialBases{n}\\
f          &\mapsto     \Conj{f(v_1)}
\end{align*}
extends to a $\SAut(F_n)$-equivariant simplicial map 
$\Phi \colon \Cayley(\SAut(F_n),S) \rightarrow \PartialBases{n}$.  When $n \geq 2$, the
group $\SAut(F_n)$ acts transitively on the vertices of $\PartialBases{n}$, so the image of
$\Phi$ contains every vertex of $\PartialBases{n}$.  Since $\Cayley(\SAut(F_n),S)$
is connected, this implies that $\PartialBases{n}$ is connected for $n \geq 2$.

Assume now that $n \geq 3$.  Our goal is to prove that $\PartialBases{n}$ is $1$-connected.  We will
prove below that the induced map
\[\Phi_{\ast} \colon \pi_1(\Cayley(\SAut(F_n),S),1) \longrightarrow \pi_1(\PartialBases{n},\Conj{v_1})\]
is surjective and has image $1$.  The desired result will then follow.

Throughout the proof, we will use the following notation.  A {\em simplicial path} in a simplicial complex
$C$ is a sequence of vertices $x_1,\ldots,x_k$ in $C$ such that for $1 \leq i < k$, the vertex
$x_i$ is either equal to $x_{i+1}$ or is joined by an edge to $x_{i+1}$.  We will denote such a path
by $x_1 - x_2 - \cdots - x_k$.
If $x_k = x_1$, then we will call this a {\em simplicial loop}.

\BeginClaims
\begin{claim}
The image of the map 
$\Phi_{\ast} \colon \pi_1(\Cayley(\SAut(F_n),S),1) \longrightarrow \pi_1(\PartialBases{n},\Conj{v_1})$
is $1$.
\end{claim}

It is well known that one can construct a $1$-connected space from the Cayley graph of a finitely-presented group by attaching discs to the orbits of the loops associated to any complete set of relations.
Let $X$ be the space obtained from $\Cayley(\SAut(F_n),S),1)$ in this way using the presentation from Theorem~\ref{theorem:gersten}.
We will show that the images in $\PartialBases{n}$ of the loops associated to these relations are contractible.  
This will imply that we can extend $\Phi$ to $X$.  Since $X$ is $1$-connected, we will be 
able to conclude that $\Phi_{\ast}$ is the zero map, as desired.

Consider one of the relations $s_1 \cdots s_k=1$ with $s_i \in S^{\pm 1}$ from Theorem \ref{theorem:gersten}.
The associated simplicial loop in $\Cayley(\SAut(F_n),S)$ and its image in $\PartialBases{n}$ are
\[1 - s_1 - s_1 s_2 - \cdots - s_1 s_2 \cdots s_k = 1\]
and
\begin{equation}
\label{eqn:tocontract}
\Conj{v_1} - \Conj{s_1(v_1)} - \Conj{s_1 s_2(v_1)} - \cdots - \Conj{s_1 s_2 \cdots s_k(v_1)} = \Conj{v_1},
\end{equation}
respectively.
Examining the relations in Theorem \ref{theorem:gersten}, we see that there exists some $x \in S_{F}$ such
that $s_i(x) = x$ for all $1 \leq i \leq k$ (this uses the fact that $n \geq 3$).  Since $\Conj{v_1}$ and
$\Conj{x}$ are either equal or joined by an edge, the vertices $\Conj{s_1 \cdots s_i(v_1)}$ and
$\Conj{s_1 \cdots s_i(x)} = \Conj{x}$ are either equal or joined by an edge for all $0 \leq i \leq k$.
We deduce that if we can show that the loops 
$\Conj{x} - \Conj{s_1 \cdots s_i(v_1)} - \Conj{s_1 \cdots s_{i+1}(v_1)} - \Conj{x}$
are contractible for all $0 \leq i < k$, then we can contract the loop in \eqref{eqn:tocontract} to $\Conj{x}$.

Using the $\SAut(F_n)$-action, we see that it is enough to show that
the loops $\Conj{x} - \Conj{v_1} - \Conj{s_i(v_1)} - \Conj{x}$ are contractible for all $1 \leq i \leq k$.
There are three cases.  In the first, $\Conj{v_1} = \Conj{s_i(v_1)}$ and the result is trivial.
In the second, the set $\{v_1, s_i(v_1), x\}$ is a partial basis for $F_n$, and again the result
is trivial.  In the third, the set $\{v_1, s_i(v_1),x\}$ is not a partial basis for $F_n$.  Looking
at the generators in $S$, we see that this can hold only if $s_i = \Mul{v_1^{e'}}{x^{e}}^{e''}$ for some
$e,e',e'' \in \{1,-1\}$.  Letting $y \in S_{F}$ be a generator distinct from $x$ and $v_1$, the
loop $\Conj{x} - \Conj{v_1} - \Conj{s_i(v_1)} - \Conj{x}$ can be contracted to $\Conj{y}$, and
we are done.

\begin{claim}
The map $\Phi_{\ast} \colon \pi_1(\Cayley(\SAut(F_n),S),1) \longrightarrow \pi_1(\PartialBases{n},\Conj{v_1})$
is surjective.
\end{claim}

Consider a loop $\Conj{w_0} - \Conj{w_1} - \cdots - \Conj{w_k}$
in $\PartialBases{n}$ with $\Conj{w_0} = \Conj{w_k} = \Conj{v_1}$.  We will show that this
loop is in the image of $\Phi_{\ast}$.
By choosing an appropriate representative $w_i$ for each conjugacy class $\Conj{w_i}$, we can assume that $w_0 = v_1$ and that
$\{w_i,w_{i+1}\}$ is a partial basis of $F_n$ for $0 \leq i < k$.  Since $\{v_1,w_1\}$ is a partial
basis for $F_n$ and $\SAut(F_n)$ acts transitively on two-element partial bases (this uses the fact
that $n \geq 3$), there exists some $\phi_1 \in \SAut(F_n,\Conj{v_1})$ such that $\phi_1(v_2) = w_1$ (in fact,
we could assume that $\phi_1$ fixes $v_1$ and not just its conjugacy class, but that is not needed).  
It follows that $w_1 = \phi_1 w_{v_2,v_1}(v_1)$.  Next, since $\{w_1,w_2\}$ is a partial basis, so
is 
\[\{(\phi_1 w_{v_2,v_1})^{-1}(w_1), (\phi_1 w_{v_2,v_1})^{-1}(w_2)\} = \{v_1,(\phi_1 w_{v_2,v_1})^{-1}(w_2)\}.\]
As before, there exists some $\phi_2 \in \SAut(F_n,\Conj{v_1})$ such that 
$\phi_2(v_2) = (\phi_1 w_{v_2,v_1})^{-1}(w_2)$.  We conclude that
$w_2 = (\phi_1 w_{v_2,v_1})(\phi_2 w_{v_2,v_1})(v_1)$.  Repeating this argument, we obtain elements
$\phi_1,\ldots,\phi_k \in \SAut(F_n,\Conj{v_1})$ such that
$w_i = (\phi_1 w_{v_2,v_1})\cdots(\phi_i w_{v_2,v_1})(v_1)$ for all $1 \leq i \leq k$. 

Set $\phi_{k+1} = ((\phi_1 w_{v_2,v_1})\cdots(\phi_i w_{v_2,v_1}))^{-1}$.  Since
\[\Conj{(\phi_1 w_{v_2,v_1})\cdots(\phi_k w_{v_2,v_1})(v_1)} = \Conj{w_k} = \Conj{v_1},\]
we have that $\phi_{k+1} \in \SAut(F_n,\Conj{v_1})$.

By Lemma \ref{lemma:stabgen}, for $1 \leq i \leq k+1$ there exists $s^i_1,\ldots,s^i_{m_i} \in S^{\pm 1}$ such
that $s^i_1 \cdots s^i_{m_i} = \phi_i$ and such that $s^i_j \in \SAut(F_n,\Conj{v_1})$ for $1 \leq j \leq m_i$.
Observe now that we have a relation
\[(s^1_1 \cdots s^1_{m_1} w_{v_2,v_1}) \cdots (s^k_1 \cdots s^k_{m_k} w_{v_2,v_1}) (s^{k+1}_1 \cdots s^{k+1}_{m_{k+1}}) = 1\]
in $\SAut(F_n)$.  The image under $\Phi_{\ast}$ of the corresponding loop in $\Cayley(\SAut(F_n),S)$ is
\[\Conj{v_1} - \Conj{s^1_1(v_1)} - \Conj{s^1_1 s^1_2(v_1)} - \cdots - \Conj{s^1_1 s^1_2 \cdots s^1_{m_1} w_{v_2,v_1}(v_1)} - \cdots.\]
Since $\Conj{s^i_j(v_1)} = \Conj{v_1}$ for all $1 \leq i \leq k+1$ and $1 \leq j \leq m_i$, after deleting
repeated vertices this path equals
\begin{align*}
\Conj{v_1} - \Conj{s^1_1 s^1_2 \cdots s^1_{m_1} w_{v_2,v_1}(v_1)} - &\Conj{(s^1_1 s^1_2 \cdots s^1_{m_1} w_{v_2,v_1})(s^2_1 s^2_2 \cdots s^2_{m_2} w_{v_2,v_1})(v_1)} - \\
&\cdots - \Conj{(s^1_1 s^1_2 \cdots s^1_{m_1} w_{v_2,v_1}) \cdots (s^k_1 s^k_2 \cdots s^k_{m_k} w_{v_2,v_1})(v_1)}.
\end{align*}
By construction, this equals $\Conj{w_0} - \Conj{w_1} - \cdots - \Conj{w_k}$, as desired.
\end{proof}

\section{The connectivity of $\PartialBases{n} / \IA_n$}
\label{section:partialbasesia}

The goal of this section is to prove Theorem \ref{theorem:basesiaconnected}.  The proof is contained in \S \ref{section:basesiaconnected}, which
is prefaced with some background information on simplicial complexes in \S \ref{section:simplicialcomplexes}.

\subsection{Simplicial complexes}
\label{section:simplicialcomplexes}

Our basic reference for simplicial complexes is \cite[Chapter 3]{Spanier}.  Let us recall the definition
of a simplicial complex given there.

\begin{definition}
A {\em simplicial complex} $X$ is a set of nonempty finite sets (called {\em simplices}) such that
if $\Delta \in X$ and $\emptyset \neq \Delta' \subset \Delta$, then $\Delta' \in X$.
The {\em dimension} of a simplex $\Delta \in X$ is $|\Delta|-1$ and is denoted $\Dim(\Delta)$.
For $k \geq 0$, the subcomplex of $X$ consisting of all simplices
of dimension at most $k$ (known as the {\em $k$-skeleton of $X$}) will be denoted $X^{(k)}$.
If $X$ and $Y$ are simplicial complexes, then a {\em simplicial map} from $X$ to $Y$ is a function
$f \colon X^{(0)} \rightarrow Y^{(0)}$ such that if $\Delta \in X$, then $f(\Delta) \in Y$.
\end{definition}

If $X$ is a simplicial complex, then we will define
the geometric realization $|X|$ of $X$ in the standard way (see \cite[Chapter 3]{Spanier}).  When
we say that $X$ has some topological property (e.g.\ simple-connectivity), we will mean that $|X|$ possesses
that property.

Next, we will need the following definitions.

\begin{definition}
Consider a simplex $\Delta$ of a simplicial complex $X$.
\begin{itemize}
\item The {\em star} of $\Delta$ (denoted $\Star_X(\Delta)$) is the subcomplex of $X$ consisting
of all $\Delta' \in X$ such that there is some $\Delta'' \in X$ with $\Delta,\Delta' \subset \Delta''$.
By convention, we will also define $\Star_X(\emptyset) = X$.
\item The {\em link} of $\Delta$ (denoted $\Link_X(\Delta)$) is the subcomplex of $\Star_X(\Delta)$
consisting of all simplices that do not intersect $\Delta$.  By convention, we will also define
$\Link_X(\emptyset) = X$.
\end{itemize}
\end{definition}

For $n \leq -1$, we will say that the empty set is both an $n$-sphere and a closed $n$-ball.  Also, if $X$
is a space then we will say that $\pi_{-1}(X)=0$ if $X$ is nonempty and that $\pi_{k}(X)=0$ for all $k \leq -2$.
With these conventions, it is true for all $n \in \Z$ that
a space $X$ satisfies $\pi_n(X)=0$ if and only if every map of an $n$-sphere into $X$ can be extended to a map
of a closed $(n+1)$-ball into $X$.

Finally, we will need the following definition.  A basic reference is \cite{RourkeSanderson}.

\begin{definition}
For $n \geq 0$, a {\em combinatorial $n$-manifold} $M$ is a nonempty simplicial complex that
satisfies the following inductive property.  If $\Delta \in M$, then $\Dim(\Delta) \leq n$.
Additionally, if $n-\Dim(\Delta)-1 \geq 0$, then $\Link_M(\Delta)$ is a combinatorial $(n-\Dim(\Delta)-1)$-manifold
homeomorphic to either an $(n-\Dim(\Delta)-1)$-sphere or a closed $(n-\Dim(\Delta)-1)$-ball.  We will denote by
$\partial M$ the subcomplex of $M$ consisting of all simplices $\Delta$ such that $\Dim(\Delta) < n$ and such that
$\Link_M(\Delta)$ is homeomorphic to a closed $(n-\Dim(\Delta)-1)$-ball.
If $\partial M = \emptyset$ then $M$  is said to be {\em closed}.  A
combinatorial $n$-manifold homeomorphic to an $n$-sphere (resp. a closed $n$-ball) will be called
a {\em combinatorial $n$-sphere} (resp. a {\em combinatorial $n$-ball}).
\end{definition}

It is well-known that if $\partial M \neq \emptyset$, then $\partial M$ is a closed combinatorial $(n-1)$-manifold
and that if $B$ is a combinatorial $n$-ball, then $\partial B$ is a combinatorial $(n-1)$-sphere.

\begin{warning}
There exist simplicial complexes that are homeomorphic to manifolds but are {\em not} combinatorial
manifolds.
\end{warning}

The following is an immediate consequence of the Zeeman's extension \cite{ZeemanSimp} of
the simplicial approximation theorem.

\begin{lemma}
\label{lemma:simpapprox}
Let $X$ be a simplicial complex and $n \geq 0$.  The following hold.
\begin{enumerate}
\item Every element of $\pi_n(X)$ is represented by a simplicial map $S \rightarrow X$,
where $S$ is a combinatorial $n$-sphere.
\item If $S$ is a combinatorial $n$-sphere and $f \colon S \rightarrow X$ is a nullhomotopic
simplicial map, then there is a combinatorial
$(n+1)$-ball $B$ with $\partial B = S$ and a simplicial map $g \colon B \rightarrow X$ such that $g|_{S} = f$.
\end{enumerate}
\end{lemma}

\subsection{The proof of Theorem \ref{theorem:basesiaconnected}}
\label{section:basesiaconnected}

We will need two lemmas.  First, some notation.
For $\{v_1,\ldots,v_k\} \subset F_n$, define
\[\Aut(F_n,v_1,\ldots,v_k) = \Set{$\phi \in \Aut(F_n)$}{$\phi(v_i) = v_i$ for $1 \leq i \leq k$}.\]
We then have the following lemma, whose proof is identical to the proof of Lemma \ref{lemma:autstabimage}.

\begin{lemma}[Images of stabilizers II]
\label{lemma:autstabimage2}
Let $\{v_1,\ldots,v_k\}$ and $\{\overline{v}_1,\ldots,\overline{v}_k\}$ be
partial bases for $F_n$ and $\Z^n$, respectively, such that
$\pi(v_i) = \overline{v}_i$ for $1 \leq i \leq k$.  Set
$\overline{V} = \Span{\overline{v}_1,\ldots,\overline{v}_k} \subset \Z^n$.  Then
the map $\Aut(F_n,v_1,\ldots,v_k) \rightarrow \Aut(\Z^n,\overline{V})$ is surjective.
\end{lemma}

\noindent
Next, we prove the following.

\begin{lemma}[Basis completion lemma]
\label{lemma:basiscompletion}
Let $\{\overline{v}_1,\ldots,\overline{v}_n\}$ be a basis for $\Z^n$ and let $\pi \colon F_n \rightarrow \Z^n$
be the abelianization map.  For some $0 \leq k \leq n$, 
let $\{v_1,\ldots,v_k\} \subset F_n$ be a partial basis for 
$F_n$ such that $\pi(v_i)=\overline{v}_i$ for $1 \leq i \leq k$.  There then 
exists $\{v_{k+1},\ldots,v_n\} \subset F_n$ such that 
$\{v_1,\ldots,v_n\}$ is a basis for $F_n$ and such that 
$\pi(v_i) = \overline{v}_i$ for $1 \leq i \leq n$.
\end{lemma}
\begin{proof}
Complete the partial basis $\{v_1,\ldots,v_k\}$ to a basis 
$\{v_1,\ldots,v_k,v_{k+1}',\ldots,v_n'\}$
for $F_n$ and set $\overline{v}_i' = \pi(v_i')$ for $k+1 \leq i \leq n$.  Define 
$\overline{V} = \Span{\overline{v}_1,\ldots,\overline{v}_k}$.  There 
then exists some $\overline{\phi} \in \Aut(\Z^n,\overline{V})$ such that
$\overline{\phi}(\overline{v}_i') = \overline{v}_i$ for $k+1 \leq i \leq n$.  By Lemma \ref{lemma:autstabimage2}, 
there exists some $\phi \in \Aut(F_n,v_1,\ldots,v_k)$
that induces $\overline{\phi}$ on $\HH_1(F_n;\Z) = \Z^n$.  The desired basis for $F_n$ 
is then 
$\{v_1,\ldots,v_k,\phi(v_{k+1}'),\ldots,\phi(v_{n}')$.
\end{proof}

\noindent
We now proceed to the proof of Theorem \ref{theorem:basesiaconnected}.

\begin{proof}[{Proof of Theorem \ref{theorem:basesiaconnected}}]
Recall that this theorem asserts that $\PartialBases{n}/\IA_n$ is $(n-2)$-connected.  Our proof will
have two steps.  In the first, we will show that $\PartialBases{n}/\IA_n$ is isomorphic to a more
concrete space $\PartialBasesZ{n}$, and in the second, we will prove that $\PartialBasesZ{n}$ is
$(n-2)$-connected.  We begin by defining $\PartialBasesZ{n}$.

\begin{definition}
Let $\PartialBasesZ{n}$ denote 
the simplicial complex whose $(k-1)$-simplices are sets $\{x_1,\ldots,x_k\} \subset \Z^n$ whose span
is a $k$-dimensional direct summand of $\Z^n$.
\end{definition}

\noindent
Now on to the proof.

\BeginSteps
\begin{step}
We have $\PartialBases{n} / \IA_n \cong \PartialBasesZ{n}$.
\end{step}

Let $\pi \colon F_n \rightarrow \Z^n$ be the projection.  There is a map 
$\psi \colon \PartialBases{n} \rightarrow \PartialBasesZ{n}$ that
takes a simplex $\{\Conj{w_1},\ldots,\Conj{w_m}\}$ of $\PartialBases{n}$ to $\{\pi(w_1),\ldots,\pi(w_m)\}$.  By Lemma \ref{lemma:basiscompletion},
the map $\psi$ is surjective.  Also, $\psi$ is invariant under the action of $\IA_n$.  It is thus enough to prove that the induced
map $\PartialBases{n}/\IA_n \rightarrow \PartialBasesZ{n}$ is injective.  To do this, it is enough to show that if $s=\{\Conj{v_1},\ldots,\Conj{v_k}\}$ and
$s'=\{\Conj{v_1'},\ldots,\Conj{v_k'}\}$ are two simplices of $\PartialBases{n}$ such that $\psi(s) = \psi(s')$, then there exists some $f \in \IA_n$ such that
$f(s) = s'$.

Reordering our simplices if necessary, we can assume that $\pi(v_i) = \pi(v_i')$ for $1 \leq i \leq k$.  Define $\overline{v}_i = \pi(v_i) \in \Z^n$.  We
can complete the partial basis $\{\overline{v}_1,\ldots,\overline{v}_k\}$ for $\Z^n$ to a basis $\{\overline{v}_1,\ldots,\overline{v}_n\}$.
Applying Lemma \ref{lemma:basiscompletion} twice, we obtain bases $\{v_1,\ldots,v_n\}$ and $\{v_1',\ldots,v_n'\}$ for $F_n$ such that
$\pi(v_i) = \pi(v_i') = \overline{v}_i$ for $1 \leq i \leq n$.  Define $f \in \Aut(F_n)$ by $f(v_i) = v_i'$.  Clearly $f(s) = s'$, and
by construction $f \in \IA_n$, as desired.

\begin{step}
The space $\PartialBasesZ{n}$ is $(n-2)$-connected.
\end{step}

This result is contained in the 1979 PhD thesis of Maazen \cite{MaazenThesis}.  Since this thesis was 
never published,
we include a proof.  Our proof is somewhat different from Maazen's proof and is modeled after proofs of 
related results due
to the second author (see \cite[Proposition 6.13]{PutmanInfinite} and \cite[Proposition 6.8]{PutmanH2}).
Fix a basis $\{e_1,\ldots,e_n\}$ for $\Z^n$, and for $0 \leq k < n$ define 
$\PartialBasesZZ{n}{k} = \Link_{\PartialBasesZ{n}}(\{e_1,\ldots,e_k\})$.  We will prove a 
more general statement, namely, that $\pi_{\ell}(\PartialBasesZZ{n}{k}) = 0$ for $0 \leq k < n$ and $-1 \leq \ell \leq n-k-2$.  

The proof will be by induction on $\ell$.  The base case $\ell=-1$ is equivalent to the observation that if $k < n$, then
$\PartialBasesZZ{n}{k}$ is nonempty.
Assume now that $0 \leq \ell \leq n-k-2$ and that
$\pi_{\ell'}(\PartialBasesZZ{n'}{k'})=0$ for all
$0 \leq k' < n'$ and $-1 \leq \ell' \leq n'-k'-2$ such that $\ell' < \ell$.  Let $S$ be a combinatorial $\ell$-sphere
and let $\phi \colon S \rightarrow \PartialBasesZZ{n}{k}$ be a simplicial map.  By Lemma \ref{lemma:simpapprox},
it is enough to show that $\phi$ may be homotoped to a constant map.

%%%%%%%%%%%%%%%%%%%
\begin{table}[ht!]
\begin{center}
\begin{tabular}{l@{\hspace{28pt}}l@{\hspace{28pt}}l}
\toprule
$\boldsymbol{t \in T}$ & $\boldsymbol{s \in S \cup S^{-1}}$ & $\boldsymbol{s t s^{-1}}$ \\
\midrule
$\Mulcomm{c}{a}{b}$	
			& $\Mul{x}{c}$      		& $\Con{x}{c}[\Con{x}{b}^{-1},\Con{x}{a}^{-1}]\Mulcomm{x}{b}{a}\Mulcomm{c}{a}{b}\Con{x}{c}^{-1}$ \\
			& $\Mul{x}{c}^{-1}$      	& $\Mulcomm{x}{a}{b}\Mulcomm{c}{a}{b}$ \\
			& $\Mul{a}{x}$      		& $\Mulcomm{c}{x}{b}\Con{c}{x}^{-1}\Mulcomm{c}{a}{b}\Con{c}{x}$ \\
			& $\Mul{a}{x}^{-1}$		& $\Con{a}{x}^{-1}\Mulcomm{c}{a}{b}\Con{c}{a}^{-1}\Con{c}{x}\Mulcomm{c}{b}{x}\Con{c}{x}^{-1}\Con{c}{a}\Con{a}{x}$ \\
			& $\Mul{a^{-1}}{x}$		& $\Mulcomm{c}{a}{b}\Con{c}{a}^{-1}\Con{c}{x}\Mulcomm{c}{b}{x}\Con{c}{x}^{-1}\Con{c}{a}$ \\
			& $\Mul{a^{-1}}{x}^{-1}$	& $\Con{a}{x}^{-1}\Mulcomm{c}{x}{b}\Con{c}{x}^{-1}\Mulcomm{c}{a}{b}\Con{c}{x}\Con{a}{x}$ \\
			& $\Mul{b}{x}$      		& $\Con{c}{x}^{-1}\Mulcomm{c}{a}{b}\Con{c}{x}\Mulcomm{c}{a}{x}$ \\
			& $\Mul{b}{x}^{-1}$      	& $\Con{b}{x}^{-1}\Con{c}{b}^{-1}\Con{c}{x}\Mulcomm{c}{x}{a}\Con{c}{x}^{-1}\Con{c}{b}\Mulcomm{c}{a}{b}\Con{b}{x}$ \\
			& $\Mul{b^{-1}}{x}$		& $\Con{c}{b}^{-1}\Con{c}{x}\Mulcomm{c}{x}{a}\Con{c}{x}^{-1}\Con{c}{b}\Mulcomm{c}{a}{b}$ \\
			& $\Mul{b^{-1}}{x}^{-1}$	& $\Con{b}{x}^{-1}\Con{c}{x}^{-1}\Mulcomm{c}{a}{b}\Con{c}{x}\Mulcomm{c}{a}{x}\Con{b}{x}$ \\
			& $\Mul{c}{x}^\epsilon$ 	& $\Con{c}{x}^\epsilon\Mulcomm{c}{a}{b}\Con{c}{x}^{-\epsilon}$ \\
			& $\Mul{a}{b}^\epsilon$			& $\Con{c}{b}^{-\epsilon}\Mulcomm{c}{a}{b}\Con{c}{b}^\epsilon$\\
			& $\Mul{a}{c}$      		& $\Con{a}{c}\Con{a}{b}\Mulcomm{a}{b}{c}\Con{a}{b}^{-1}\Con{a}{c}^{-1}\Con{a}{b}[\Con{c}{a}^{-1},\Con{c}{b}^{-1}]\Mulcomm{c}{a}{b}\Con{c}{b}^{-1}$ \\
			& $\Mul{a}{c}^{-1}$      	& $\Con{c}{b}\Mulcomm{c}{a}{b}\Con{c}{a}^{-1}\Con{a}{c}\Mulcomm{a}{b}{c}\Con{a}{b}^{-1}\Con{a}{c}^{-1}\Con{c}{a}$ \\
			& $\Mul{a^{-1}}{c}$		& $\Con{a}{c}\Con{c}{b}\Mulcomm{c}{a}{b}\Con{c}{a}^{-1}\Con{a}{c}\Mulcomm{a}{b}{c}\Con{a}{b}^{-1}\Con{a}{c}^{-1}\Con{c}{a}\Con{a}{c}^{-1}$ \\
			& $\Mul{a^{-1}}{c}^{-1}$	& $\Con{a}{b}\Mulcomm{a}{b}{c}\Con{a}{b}^{-1}\Con{a}{c}^{-1}\Con{a}{b}[\Con{c}{a}^{-1},\Con{c}{b}^{-1}]\Mulcomm{c}{a}{b}\Con{c}{b}^{-1}\Con{a}{c}$ \\
			& $\Mul{b}{a}^\epsilon$			& $\Con{c}{a}^{-\epsilon}\Mulcomm{c}{a}{b}\Con{c}{a}^\epsilon$\\
			& $\Mul{b}{c}$      		& $\Con{c}{a}\Mulcomm{c}{a}{b}[\Con{c}{a}^{-1},\Con{c}{b}^{-1}]\Con{b}{a}^{-1}\Con{b}{c}\Con{b}{a}\Mulcomm{b}{c}{a}\Con{b}{a}^{-1}\Con{b}{c}^{-1}$ \\
			& $\Mul{b}{c}^{-1}$     	& $\Con{c}{b}^{-1}\Con{b}{c}\Con{b}{a}\Mulcomm{b}{c}{a}\Con{b}{c}^{-1}\Con{c}{b}\Mulcomm{c}{a}{b}\Con{c}{a}^{-1}$ \\
			& $\Mul{b^{-1}}{c}$		& $\Con{b}{c}\Con{c}{b}^{-1}\Con{b}{c}\Con{b}{a}\Mulcomm{b}{c}{a}\Con{b}{c}^{-1}\Con{c}{b}\Mulcomm{c}{a}{b}\Con{c}{a}^{-1}\Con{b}{c}^{-1}$ \\
			& $\Mul{b^{-1}}{c}^{-1}$	& $\Con{b}{c}^{-1}\Con{c}{a}\Mulcomm{c}{a}{b}[\Con{c}{a}^{-1},\Con{c}{b}^{-1}]\Con{b}{a}^{-1}\Con{b}{c}\Con{b}{a}\Mulcomm{b}{c}{a}\Con{b}{a}^{-1}$ \\
			& $\Mul{c}{a}^\epsilon$ 	& $\Con{c}{a}^\epsilon\Mulcomm{c}{a}{b}\Con{c}{a}^{-\epsilon}$ \\
			& $\Mul{c}{b}^\epsilon$ 	& $\Con{c}{b}^\epsilon\Mulcomm{c}{a}{b}\Con{c}{b}^{-\epsilon}$ \\

\midrule
$\Con{c}{a}$	
		& $\Mul{x}{c}$            	& $\Con{c}{a}\Con{x}{c}\Con{x}{a}\Mulcomm{x}{c}{a}\Con{x}{a}^{-1}\Con{x}{c}^{-1}$ \\
		& $\Mul{x}{c}^{-1}$            	& $\Con{c}{a}\Con{x}{a}\Mulcomm{x}{a}{c}\Con{x}{a}^{-1}$ \\
		& $\Mul{x^{-1}}{c}$      	& $\Con{x}{c}\Con{c}{a}\Con{x}{a}\Mulcomm{x}{a}{c}\Con{x}{a}^{-1}\Con{x}{c}^{-1}$ \\
		& $\Mul{x^{-1}}{c}^{-1}$      	& $\Con{x}{c}^{-1}\Con{c}{a}\Con{x}{c}\Con{x}{a}\Mulcomm{x}{c}{a}\Con{x}{a}^{-1}$ \\
		& $\Mul{a^\epsilon}{x}^\delta$ 	& $(\Con{c}{a}^\epsilon\Con{c}{x}^\delta)^\epsilon$ \\
		& $\Mul{c}{x}$            	& $\Con{c}{a}\Con{c}{x}\Mulcomm{c}{a}{x}\Con{c}{x}^{-1}$ \\
		& $\Mul{c}{x}^{-1}$ 		& $\Con{c}{a}\Mulcomm{c}{x}{a}$\\
		& $\Mul{c^{-1}}{x}$      	& $\Con{c}{x}\Con{c}{a}\Mulcomm{c}{x}{a}\Con{c}{x}^{-1}$ \\
		& $\Mul{c^{-1}}{x}^{-1}$  	& $\Con{c}{x}^{-1}\Con{c}{a}\Con{c}{x}\Mulcomm{c}{a}{x}$ \\
		& $\Mul{a^\epsilon}{c}^\delta$ 	& $(\Con{a}{c}^\delta\Con{c}{a}^\epsilon)^\epsilon$ \\
\bottomrule
\end{tabular}
\caption{Conjugating $T$ by $S \cup S^{-1}$.
Here $a$, $b$, $c$ and $x$ are distinct elements of the generating set for $F_n$ 
and $\epsilon,\delta=\pm 1$.
If a particular choice of $s$ does not appear in the table, then $sts^{-1}=s^{-1}ts=t$.}
\label{table:AutactiononIAn}
\end{center}
\end{table}
%%%%%%%%%%%%%%%%%%%

We begin with some notation.  Consider $w \in \Z^n$.  Write $w = \sum_{i=1}^n c_i e_i$ with $c_i \in \Z$, and define
$\Rank(w) = |c_n|$.  Now set
\[R = \Max \Set{$\Rank(\phi(x))$}{$x \in S^{(0)}$}.\]
If $R = 0$, then $\phi(S) \subset \Star_{\PartialBasesZZ{n}{k}}(\{e_n\})$,
and hence the map $\phi$ can be homotoped to the constant map $e_n$. 
Assume, therefore, that $R>0$.  Let $\Delta'$ be a simplex of $S$ such that
$\Rank(\phi(x))=R$ for all vertices $x$ of $\Delta'$.  Choose
$\Delta'$ so that $m := \Dim(\Delta')$ is maximal, which implies that $\Rank(\phi(x)) < R$
for all vertices $x$ of $\Link_S(\Delta')$.  Now, $\Link_S(\Delta')$ is a combinatorial $(\ell-m-1)$-sphere and
$\phi(\Link_S(\Delta'))$ is contained in
\[\Link_{\PartialBasesZZ{n}{k}}(\phi(\Delta')) \cong \PartialBasesZZ{n}{k+m'}\]
for some $m' \leq m$ (it may be less than $m$ if $\phi|_{\Delta'}$ is not injective).  The inductive hypothesis together
with Lemma \ref{lemma:simpapprox} therefore tells us that there a combinatorial $(\ell-m)$-ball $B$ with
$\partial B = \Link_S(\Delta')$ and a simplicial map $f \colon B \rightarrow \Link_{\PartialBasesZZ{n}{k}}(\phi(\Delta'))$ such that
$f|_{\partial B} = \phi|_{\Link_S(\Delta')}$.

Our goal now is to adjust $f$ so that $\Rank(\phi(x)) < R$ for all $x \in B^{(0)}$.
Let $v \in \Z^n$ be a vector corresponding to a vertex of $\phi(\Delta')$.  Observe that the $e_n$-coordinate of $v$ is
$\pm R$.
We define a map $f' \colon B \rightarrow \Link_{\PartialBasesZZ{n}{k}}(\phi(\Delta'))$ in the following way.  Consider
$x \in B^{(0)}$, and let $v_x = f(x) \in \Z^n$.
By the division algorithm, there exists some $q_x \in \Z$
such that $\Rank(v_x + q_x v) < R$.  Moreover, by the maximality of $m$ we can choose $q_x$
such that $q_x = 0$ if $x \in (\partial B)^{(0)}$.
Define $f'(x) = v_x + q_x v$.  It is clear that the map
$f'$ extends to a map $f' \colon B \rightarrow \Link_{\PartialBasesZZ{n}{k}}(\phi(\Delta'))$.  Additionally,
$f'|_{\partial B} = f|_{\partial B} = \phi|_{\Link_S(\Delta')}$.
We conclude that we can homotope $\phi$ so as to replace $\phi|_{\Star_S(\Delta')}$ with $f'$.  Since
$\Rank(f'(x)) < R$ for all $x \in B$, we have removed $\Delta'$ from $S$ without introducing any
vertices whose images have rank greater than or equal to $R$.  Continuing in this manner allows
us to simplify $\phi$ until $R=0$, and we are done.
\end{proof}

\appendix
\section{Appendix : Derivation of Theorem \ref{theorem:magnusfinite} from Theorem \ref{theorem:magnus}}
\label{appendix:magnus}

Fixing a basis $\{v_1,\ldots,v_n\}$ for $F_n$, let $T$ be the purported generating set 
for $\IA_n$ from Theorem \ref{theorem:magnusfinite} and
let $\Gamma < \Aut(F_n)$ be the subgroup generated by $T$.  Fix some $f \in \Aut(F_n)$.
By Theorem \ref{theorem:magnus}, to show that $\Gamma = \Aut(F_n)$, it is enough to show that
$f \Con{v_1}{v_2} f^{-1} \in \Gamma$.  Recall that $\SAut(F_n)$ consists
of all elements of $\Aut(F_n)$ whose images in $\Aut(\Z^n)$ have determinant $1$.  Letting 
$I_1 \in \Aut(F_n)$ be the automorphism that takes $v_1$ to $v_1^{-1}$ and fixes $v_i$ for $i > 1$,
we can write $f = g \cdot I_1^{k}$ for some $g \in \SAut(F_n)$ and some $k \in \Z$.  We then have
$$f \Con{v_1}{v_2} f^{-1} = g \cdot I_1^k \Con{v_1}{v_2} I_1^{-k} \cdot g^{-1} = g \Con{v_1}{v_2} g^{-1}.$$
Since $\Con{v_1}{v_2} \in \Gamma$, to prove that $g \Con{v_1}{v_2} g^{-1} \in \Gamma$ it
is enough to prove that $\Gamma$ is a normal subgroup of $\SAut(F_n)$.  
Define
$$S = \Set{$\Mul{v_i}{v_j}$}{$1 \leq i,j \leq n$, $i \neq j$}.$$
The group $\SAut(F_n)$ is generated by $S$ (see Theorem \ref{theorem:gersten}), so it is enough
to show that for $s \in S \cup S^{-1}$ and $t \in T$, the automorphism $s t s^{-1}$ can
be written as a word in $T \cup T^{-1}$.  The various cases of this are contained in
Table \ref{table:AutactiononIAn}.

\noindent
\begin{tabular*}{\linewidth}[t]{@{}p{\widthof{Department of Mathematics 253-37}+1in}@{}p{\linewidth - \widthof{Department of Mathematics 253-37}-1in}@{}}
{\raggedright
Matthew Day\\
Department of Mathematics 253-37\\
California Institute of Technology\\
Pasadena, CA 91125\\
E-mail: {\tt mattday@caltech.edu}}
&
{\raggedright
Andrew Putman\\
Department of Mathematics\\
Rice University, MS 136 \\
6100 Main St.\\
Houston, TX 77005\\ 
E-mail: {\tt andyp@rice.edu}}
\end{tabular*}

\end{document}